\documentclass[10pt,a4paper]{article}
\usepackage{amsfonts}
\usepackage{latexsym}
\usepackage{cite}
\usepackage{amsmath,amsfonts,latexsym,amssymb}
\usepackage[mathscr]{eucal}
\usepackage{cases,color}
\usepackage{amsthm}

\usepackage[bf,small]{caption2}
\usepackage{float}
\usepackage{graphicx}
\usepackage{amsmath}
\usepackage{amssymb}
\usepackage[all]{xy}

\newtheorem{theorem}{theorem}[section]
\newtheorem{thm}[theorem]{Theorem}
\newtheorem{lem}[theorem]{Lemma}
\newtheorem{prob}[theorem]{Problem}

\newtheorem{cor}[theorem]{Corollary}

\newtheorem{case}[theorem]{Case}

\newtheorem{exmp}[theorem]{Example}

\newtheorem{rmk}[theorem]{Remark}
\newtheorem{nota}[theorem]{Notation}

\begin{document}

\title{\textbf{Alexander polynomial, Dijkgraaf-Witten invariant, and Seifert fibred surgery}}
\author{\Large Haimiao Chen}

\date{}
\maketitle

\begin{abstract}
  We apply Dijkgraaf-Witten invariants over semidirect products of abelian groups to show that, if the $k/\ell$-surgery along a knot $K$ results in a small Seifert 3-manifold with multiplicities $a_1,a_2,a_3$, then many constraints on $k,a_1,a_2,a_3$ can be read off from the Alexander polynomial of $K$.

\end{abstract}

\section{Introduction}

For a knot $K\subset S^3$, let $N(K)$ denote a tubular neighborhood. Given $k/\ell\in\mathbb{Q}$,
the {\it $k/\ell$-surgery} on $K$ is the operation of constructing the closed manifold
$$K(k/\ell):=(S^3-N(K))\cup_{f(k/\ell)}D^2\times S^1,$$
using a gluing map $f(k/\ell):S^1\times S^1=\partial(D^2\times S^1)\to\partial N(K)=S^1\times S^1$ which sends $S^1\times 1$ to a curve of slope $k/\ell$.

It was known to Thurston \cite{Th80} (Theorem 5.8.2) in 1980 that there are finitely many exceptional surgeries for each hyperbolic $K$.
Agol \cite{Ag98} (1998) and M. Lackenby \cite{La00} (2000) proved that a hyperbolic knot admits at most 12 exceptional surgeries. Later, the number is reduced to 10 by Lackenby and Meyerhoff \cite{LM13} (2013).
By Perelman's work, the result of an exceptional surgery is either reducible, or toroidal, or Seifert fibred.
A famous conjecture of Gordon \cite{Go09} asserts that if $K(k/\ell)$ is small Seifert fibred, then $\ell=\pm 1$, (see Problem 1.77 in \cite{Ki97} for related topics.)
This was confirmed for alternating knots in \cite{Ic08} (2008).
Works of Wu etc. \cite{Wu11, Wu11-2, Wu13, IM16} completely determined exceptional surgeries of arborescent knots, in particular, verifying Gordon's conjecture for these knots.
%Besides, the work \cite{IM18} is a recent proceeding.

Kadokami \cite{Ka07} (2007) obtained constraints on $k$ for $K(k/\ell)$ to be small Seifert, in terms of $\Delta_K$, the Alexander polynomial of $K$.
Using knot Floer homology, Wu \cite{WuZ12} (2012) gave obstructions for Seifert fibred surgeries; in particular, restrictions on the coefficients of $\Delta_K$ were induced.

We also aim to utilize Alexander polynomial to give necessary conditions, as far as possible.

Given a positive integer $m$ and a polynomial $\alpha\in\mathbb{Z}[t]$, let
\begin{align*}
\Omega_m(\alpha)&=\prod\limits_{s\in(\mathbb{Z}/m\mathbb{Z})^\times}|\alpha(\zeta_m^s)|, \qquad \zeta_m=e^{2\pi\sqrt{-1}/m},  \\
\widetilde{\Omega}_m(\alpha)&=\prod\limits_{0\ne s\in\mathbb{Z}/m\mathbb{Z}}|\alpha(\zeta_m^s)|=\prod\limits_{m'\mid m}\Omega_{m'}(\alpha); \\
[m]_t&=1+t+\cdots+t^{m-1}.
\end{align*}
%Let $[m]_t=1+t+\cdots+t^{m-1}$.

The main result of this paper is
\begin{thm} \label{thm:main}
Suppose $k\ge 2$, and $K(k/\ell)$ is small Seifert fibred 3-manifold with multiplicities $a_1,a_2,a_3\ge 2$. Let $e_1=(a_2,a_3), e_2=(a_1,a_3), e_3=(a_1,a_2)$.
\begin{enumerate}
  \item[\rm(i)] We have the greatest common divisor $(a_1,a_2,a_3)=1$.
  \item[\rm(ii)] If $\widetilde{\Omega}_k(\Delta_K)=1$, then $e_1=e_2=e_3=1$.
  \item[\rm(iii)] If $e_j>1$ for some $j$, then $\Delta_K(t)\equiv a_j\alpha(t)\pmod{[e_j]_t}$ for some $\alpha\in\mathbb{Z}[t]$ such
         that $\deg\alpha\le e_j-2$ and $\widetilde{\Omega}_{e_j}(\alpha)=1$.
  \item[\rm(iv)] If $\widetilde{\Omega}_k(\Delta_K)=0$, then $e_j>1$ for at least two $j$'s,
        $$\frac{[e_1e_2e_3]_t}{[e_1]_t[e_2]_t[e_3]_t}\ \Big|\ \Delta_K(t),$$
        and $d\mid e_1e_2e_3$ for each divisor $d$ of $k$ with $\Omega_d(\Delta_K)=0$.
  \item[\rm(v)] If $\widetilde{\Omega}_k(\Delta_K)>1$, then $e_j>1$ for exactly one $j$, and $d\mid e_j$ for each divisor $d$ of $k$ with $\Omega_d(\Delta_K)>1$.
\end{enumerate}
\end{thm}

\begin{rmk}\label{rmk:k-determine-e}
\rm In the case $\widetilde{\Omega}_k(\Delta_K)>1$, when $k$ is given, $e_j$ is determined by that it is the least common multiple of all $d$'s with $\Omega_d(\Delta_K)>1$, and $a_j$ is determined via $a_j^{\phi(e_j)}=\Omega_{e_j}(\Delta_K)$ which is implied by (iii).

Furthermore, there is an upper bound $e_j\le\deg\Delta_K+1$: otherwise $\Delta_K=a_j\alpha$, but, as well-known, $\Delta_K(1)\in\{\pm 1\}$.
\end{rmk}

This work is motivated by realizing the fact that Alexander polynomial, which is closely related to Reidemeister torsion, is in turn closely related to Dijkgraaf-Witten invariant over a semiproduct $A\rtimes(\mathbb{Z}/n\mathbb{Z})$, with $A$ abelian. When $K(k/\ell)=S^2(b_1/a_1,b_2/a_2,b_3/a_3)=M$, the DW invariant of $M$, which enumerates homomorphisms $\pi_1(M)\to A\rtimes(\mathbb{Z}/n\mathbb{Z})$, can be computed in two ways: one is expressed in terms of the Alexander module of $K$ (so involving $\Delta_K$) and $k$, and the other is given by applying the formula derived by the author \cite{Ch12} for general Seifert 3-manifolds, expressed in terms of $k,a_1,a_2,a_3$, etc. In this way, we obtain ``strange" equalities, which have many nontrivial implications.
We re-derive the results of Kadokami \cite{Ka07} (Theorem 1.1--1.4 there), and actually enhance them.

For a more broad background, in recent years, many deep connections between 3-manifolds and finite groups have been uncovered (see \cite{LR10,WZ17,WZ19} and the references therein). Remarkably, Ueki \cite{Ue18} (2018) showed that $\Delta_K$ is determined by quotients of $\pi_1(S^3-N(K))$ onto finite metabelian groups. %It should be meaningful to utilize Alexander polynomial and finite groups as far as possible.
In the context of Seifert fibred surgery, with the metabelian groups in the form $A\rtimes(\mathbb{Z}/n\mathbb{Z})$, we reveal some fine algebraic structures, enriching the topic of studying 3-manifolds via finite groups.

The content is organized as follows.
In Section 2 we clarify a connection between DW invariant of $K(k/\ell)$ and $\Delta_K$. In Section 3 we give a formula for the DW invariant of a small Seifert 3-manifold, with the proof delayed to Section 6.
In Section 4 we prove Theorem \ref{thm:main}. In Section 5, we give some examples to illustrate the utility of Theorem \ref{thm:main}.

\begin{nota}
\rm For $k\in\mathbb{Z}$, let $\Upsilon(k)$ denote the set of prime divisors of $k$. Let $\mathbb{Z}_{\ge k}=\{a\in\mathbb{Z}\colon a\ge k\}$.

Given a prime number $p$. For $h\ge 1$, let $\mathbb{F}_{p^h}$ denote the field with $p^h$ elements. For $0\ne a\in\mathbb{Z}$, let $\|a\|_p=\sup\{s\in\mathbb{Z}_{\ge 0}\colon p^s\mid a\}$. %set $\|0\|_p=-\infty$.

For a finite set $X$, let $\#X$ denote its cardinality.

Given $n\in\mathbb{Z}_{\ge 1}$. Let $\Psi_n(t)$ denote the $n$-th cyclotomic polynomial; let $\phi(n)=\#\{s\colon 1\le s<n,\ (s,n)=1\}=\deg\Psi_n$. %Let $\zeta_n=\exp(2\pi i/n)$.
Let $\mathbb{Z}_n=\mathbb{Z}/n\mathbb{Z}$, regarded as a quotient ring of $\mathbb{Z}$; when $n=p$ is a prime, $\mathbb{Z}_p$ is identified with $\mathbb{F}_p$.
For $a\in\mathbb{Z}$, denote its image under the quotient map $\mathbb{Z}\twoheadrightarrow\mathbb{Z}_n$ also by $a$.

For elements $x, y$ of some group, let $y.x=yxy^{-1}$, let
${\rm Cen}(x)$ denote the centralizer of $x$, and let ${\rm Con}(x)$ denote the conjugacy class containing $x$.

For a condition $\mathfrak{c}$, set $\delta_{\mathfrak{c}}=1$ (resp. $\delta_{\mathfrak{c}}=0$) if $\mathfrak{c}$ holds (resp. does not hold);
for instance, $\delta_{2<3}=1$, $\delta_{3|5}=0$.
\end{nota}

\section{Dijkgraaf-Witten invariant over a semidirect product and Alexander polynomial}

Suppose $A$ is a finite abelian group, $f\in{\rm Aut}(A)$ has order $n$, and suppose ${\rm id}-f^v$ (as an endomorphism of $A$) is invertible for any $v\not\equiv 0\pmod{n}$.
Let
\begin{align*}
G=A\rtimes_f\mathbb{Z}_n,
\end{align*}
i.e., the semidirect product determined by the homomorphism $\mathbb{Z}_n\to{\rm Aut}(A)$ sending $v$ to $f^v$.
Write an element of $G$ as $u\beta^v$, with $u\in A$ and $v\in\mathbb{Z}_n$. Then
\begin{align}
u_1\beta^{v_1}u_2\beta^{v_2}&=(u_1+f^{v_1}u_2)\beta^{v_1+v_2}, \nonumber \\
(u_1\beta^{v_1}).(u_2\beta^{v_2})&=(({\rm id}-f^{v_2})u_1+f^{v_1}u_2)\beta^{v_2},  \label{eq:metabelian-2} \\
(u_1\beta^{v_1})^{-1}.(u_2\beta^{v_2})&=((f^{v_2-v_1}-f^{-v_1})u_1+f^{-v_1}u_2)\beta^{v_2}.  \label{eq:metabelian-3}
\end{align}
So the commutator subgroup $G'=A$, and the abelianization $G^{\rm ab}\cong\mathbb{Z}_n.$

Note that if $v\ne 0$, then for each $s\in\mathbb{Z}_{\ge 1}$,
\begin{align}
(u\beta^v)^s=([s]_{f^v}u)\beta^{sv}=(({\rm id}-f^v)^{-1}({\rm id}-f^{sv})u)\beta^{sv},  \label{eq:power}
\end{align}
and ${\rm Cen}(u\beta^v)=\big\{u'\beta^{v'}\colon({\rm id}-f^v)u'=\big({\rm id}-f^{v'}\big)u\big\}$.

Take a projection diagram for $K$, with directed arcs $\tilde{x}_1,\ldots,\tilde{x}_b$. There is an associated Wirtinger presentation
$$\pi(K):=\pi_1(S^3-N(K))=\langle \tilde{x}_1,\ldots,\tilde{x}_b\mid \tilde{r}_1,\ldots,\tilde{r}_{b-1}\rangle,$$
with $\tilde{r}_i\in F_b$, the free group generated by $\tilde{x}_1,\ldots,\tilde{x}_b$. Let $x_j\in\pi(K)$ denote the element represented by $\tilde{x}_j$.
Recall one of the definitions of Alexander polynomial:
Consider the ring homomorphism
$$\Theta:\mathbb{Z}[F_b]\to\mathbb{Z}[\pi(K)]\to\mathbb{Z}[t^{\pm 1}]$$
sending $\tilde{x}_j^{\pm 1}$ to $t^{\pm 1}$, and put
$$\mathcal{M}=\left(\Theta(\partial\tilde{r}_i/\partial\tilde{x}_j)\right)_{(b-1)\times b},$$
where $\partial\tilde{r}_i/\partial\tilde{x}_j$ is the Fox derivative of $\tilde{r}_i$ with respect to $\tilde{x}_j$.
Let $\mathcal{M}'$ be the matrix obtained from deleting the last column of $\mathcal{M}$. Then $\Delta_K(t)=\pm t^s\det(\mathcal{M}')\in\mathbb{Z}[t^{\pm1}]$, where the sign and $s$ are chosen to ensure that
$\Delta_K(t)$ is a polynomial in $t$ with $\Delta_K(0)\ne 0$ and $\Delta_K(1)=1$ (as a convention taken in this paper).

\begin{figure}[h]
  \centering
  \includegraphics[width=7cm]{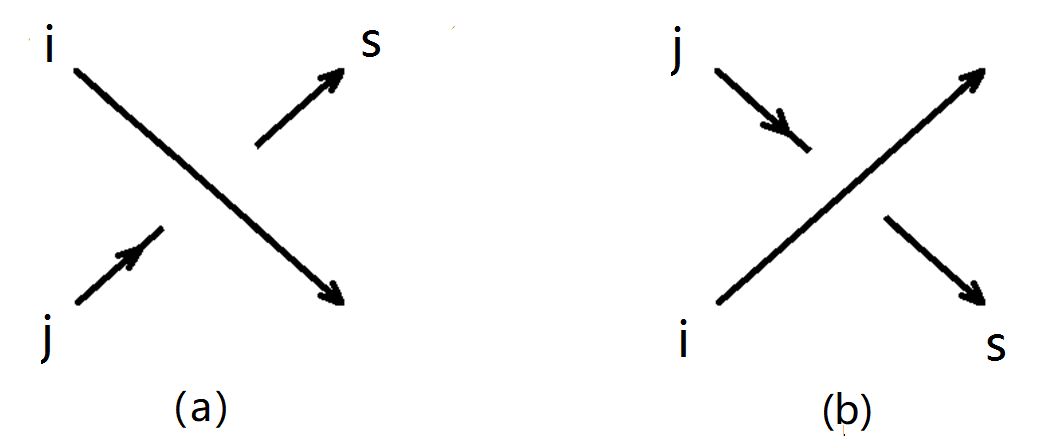}\\
  \caption{(a) a positive crossing; (b) a negative crossing.} \label{fig:pm}
\end{figure}

Given a homomorphism $\varrho:\pi(K)\to G$, since the $x_j$'s are conjugate to each other, we can assume $\varrho(x_j)=u_j\beta^{v}$, $j=1,\ldots,b$.
By (\ref{eq:metabelian-2}), (\ref{eq:metabelian-3}), if the directed arcs $\tilde{x}_i,\tilde{x}_j,\tilde{x}_s$ form a positive (resp. negative) crossing as in Figure \ref{fig:pm},
then
$$u_s=(1-f^{v})u_i+f^{v}u_j, \qquad (\text{resp.}\ \ u_s=(1-f^{-v})u_i+f^{-v}u_j).$$
So $\mathcal{M}'_{v}(u_1-u_b,\ldots,u_{b-1}-u_b)^T=0$, where $\mathcal{M}'_{v}$ is obtained from $\mathcal{M}'$ by replacing $t$ with $f^{v}$, and $T$ stands for transpose.

Conversely, with $u_b\in A$ and $v\in\mathbb{Z}_n$ fixed, each $\vec{u}\in A^{b-1}$ satisfying $\mathcal{M}'_{v}\vec{u}=0$ gives rise to a unique homomorphism $\varrho: \pi(K)\to G$.

Choose $x_b$ as the meridian for surgery. Since the element represented by the longitude lies in the second derived group of $\pi(K)$, we have that $\varrho$ extends to $\pi_1(K(k/\ell))\to G$ if and only if $(u_b\beta^v)^k=1$, which is equivalent to $v=0, ku_b=0$ or, according to (\ref{eq:power}), $v\ne 0=kv$.
Thus
\begin{align*}
\#\hom(\pi_1(K(k/\ell)),G)=\#\{u\in A\colon ku=0\}+\#A\sum\limits_{0\ne v\in\mathbb{Z}_n\atop kv=0}\#\ker\mathcal{M}'_v.
\end{align*}

We mainly focus on two cases for $G=A\rtimes_f\mathbb{Z}_n$.
\begin{case}\label{case:1}
\rm Suppose $p$ is a prime, $h\in\mathbb{Z}_{\ge 1}$, $n\in\mathbb{Z}_{\ge 2}$, such that $p\nmid n$ and $n\mid p^h-1$. Let $A=\mathbb{F}_{p^h}$, which is isomorphic to $\mathbb{Z}_p^h$ as an abelian group, and let $f\in{\rm Aut}(A)$ be given by multiplication by a primitive $n$-th root of unity $\xi_n\in\mathbb{F}_{p^h}^\times$.
Note that $\mathcal{M}'_{v}\vec{u}=0$ is a linear system of equations over $\mathbb{F}_{p^h}$; let $\omega(v)$ denote the dimension of the solution space. %Clearly, $\omega(v)>0$ if and only if $\Delta_K(\xi_n^{v})=0\in\mathbb{F}_{p^h}$.
Then
\begin{align}
\#\hom(\pi_1(K(k/\ell)),G)=p^{h\delta_{p\mid k}}+p^h\sum\limits_{j=1}^{(n,k)-1}p^{h\omega(jn/(n,k))}. \label{eq:equality-1'}
\end{align}
\end{case}

%The second case requires a preliminary.

\begin{case}\label{case:2}
\rm  Given $g\in\mathbb{Z}_{\ge 2}$ and a prime $p\nmid g$, choose a set $\Xi_p(g)$ of coset representatives for the subgroup $\langle p\rangle\le\mathbb{Z}_g^\times$; denote $\#\langle p\rangle$ by $o_p(g)$ (or simply $o(g)$). According to \cite{LN83} Theorem 2.45 and Theorem 2.47, over the algebraic closure $\overline{\mathbb{F}_p}$ of $\mathbb{F}_p$, $\Psi_g(t)=\prod_{i\in\mathbb{Z}_g^\times}(t-\xi_g^i),$
where $\xi_g$ is a primitive $g$-th root of unity in $\overline{\mathbb{F}_p}$, and over $\mathbb{F}_p$,
%$\Psi_g$ is decomposed into $\phi(g)/o_p(g)$ monic irreducible polynomials of degree $o_p(g)$:
$$\Psi_g(t)=\prod\limits_{s\in\Xi_p(g)}\psi_{g,s}, \qquad \text{with} \ \ \ \psi_{g,s}=\prod\limits_{i=1}^{o_p(g)}\big(t-\xi_g^{sp^i}\big)\in\mathbb{Z}_p[t].$$
These $\psi_{g,s}$'s are irreducible and pairwise coprime.
For $m\ge 2$, by Hensel's Lemma, each $\psi_{g,s}$ has a unique lift to $\mathbb{Z}_{p^m}[t]$ which we denote by the same notation, such that the above factorization holds in $\mathbb{Z}_{p^m}[t]$.

Let $A=A^{(m)}=\mathbb{Z}_{p^m}[t]/(\psi_{g,s}(t))$ which is isomorphic to $\mathbb{Z}_{p^m}^h$ as an abelian group, where $h=o_p(g)$;
let $f$ be given by multiplication by $t$. Then $t^v-1$ is invertible in $A$ for each $0\ne v\in\mathbb{Z}_g$, as $(t^v-1,\psi_{g,s}(t))=1$ in $\mathbb{Z}_{p^m}[t]$. We have
\begin{align}
\#\hom(\pi_1(K(k/\ell)),G)=(p^m,k)^{h}+p^{mh}\sum\limits_{j=1}^{(g,k)-1}p^{h\varpi_m(jg/(g,k))}. \label{eq:equality-2'}
\end{align}
Here $\varpi_m(v)$ is defined as follows. Over $A^{(1)}$ which is a field, $\mathcal{M}'_v$ can be converted into a diagonal matrix
${\rm diag}(1,\ldots,1,0,\ldots,0)$ (the number of zero's denoted by $b_0$) through elementary row- and column-transformations.
It is easy to see that an element of $A^{(m)}$ is invertible if and only if its image under the canonical map
$A^{(m)}\twoheadrightarrow A^{(1)}$ is invertible. Hence, over $A^{(m)}$, $\mathcal{M}'_v$ can be converted into ${\rm diag}(1,\ldots,1,p^{i(1)},\ldots,p^{i(b_0)})$ for some $i(1),\ldots,i(b_0)\in\{1,\ldots,m\}$. Put $\varpi_m(v)=i(1)+\cdots+i(b_0)$. Indeed $\#\ker\mathcal{M}'_v=p^{h\varpi_m(v)}$.
\end{case}

\section{The Dijkgraaf-Witten invariant over $\mathbb{Z}_{p^m}^h\rtimes_f\mathbb{Z}_n$ of a small Seifert 3-manifold}

According to \cite{Sa99} Section 1.6, a Seifert 3-manifold of genus 0 can be constructed by gluing $n$ copies of $D^2\times S^1$ onto $M'_n:=(S^2-\sqcup_nD^2)\times S^1$, via diffeomorphisms
$f_j:S^1\times S^1=\partial(D^2\times S^1)\to\partial M'_n=\sqcup_n(S^1\times S^1)$, $j=1,\ldots,n$.
When $f_j$ represents the mapping class
$$\left(\begin{array}{cc} a_j & b_j \\ a'_j & b'_j \end{array}\right)\in{\rm SL}(2,\mathbb{Z}), \qquad a_j\ge 2,$$
the resulting manifold is denoted by $S^2(b_1/a_1,\ldots,b_n/a_n)$, and the image of $\{0\}\times S^1$ in the $j$-th $D^2\times S^1$ under the gluing is called the $j$-th {\it singular fiber}, with {\it multiplicity} $a_j$.
In particular, $S^2(b_1/a_1,b_2/a_2,b_3/a_3)$ is called a {\it small Seifert 3-manifold}.

\begin{thm} \label{thm:DW-formula}
Let $z=\max\{\|a_1\|_p,\|a_2\|_p,\|a_3\|_p\}$, $e_j=(a_{j-1},a_{j+1})$, (with cyclical notations adopted), and
\begin{align*}
\sum\limits_{j=1}^3\frac{b_j}{a_j}&=\frac{\mathfrak{n}}{[a_1,a_2,a_3]}, \qquad \mathfrak{d}=([a_1,a_2,a_3],n), \qquad  \mathfrak{e}=(n/\mathfrak{d},\mathfrak{n}), \\
\mathfrak{k}_2&=\mathfrak{k}_2(n)=2-\sum\limits_{j=1}^3(e_{j},n)+\frac{1}{\mathfrak{d}}\prod\limits_{j=1}^3(a_j,n), \\
\mathfrak{k}_1&=\mathfrak{k}_1(n,p^m,h)
=\sum\limits_{j=1}^3((e_j,n)-1)(a_j,p^m)^h+\frac{\mathfrak{e}-1}{\mathfrak{d}}\prod\limits_{j=1}^3(a_j,n), \\
\mathfrak{k}_0&=\mathfrak{k}_0(p^m,h)=\begin{cases} p^{h\min\{m,\|\mathfrak{n}\|_p\}}, &z=0, \\
p^{-h\min\{m,z\}}\prod_{j=1}^3(a_j,p^m)^h, &z>0.\end{cases}
\end{align*}
Then
\begin{align*}
\#\hom\big(\pi_1\big(S^2(b_1/a_1,b_2/a_2,b_3/a_3)\big),\mathbb{Z}_{p^m}^h\rtimes_f\mathbb{Z}_n\big)=\mathfrak{k}_2p^{2mh}+\mathfrak{k}_1p^{mh}+\mathfrak{k}_0.
\end{align*}
\end{thm}

This formula should have independent value. %We delay the proof to the final section.

%\begin{rmk} \label{rmk:k2}
%\rm Let $q$ be a prime. Suppose $\{\|a_1\|_q,\|a_2\|_q,\|a_3\|_q\}=\{c_1,c_2,c_3\}$ with $c_1\le c_2\le c_3$. Then $\mathfrak{k}_2(q^c)>0$ only when $c_1>0$, i.e. $q\mid(a_1,a_2,a_3)$, in which case
%$$\mathfrak{k}_2(q^c)=(q^{\min\{c,c_1\}}-1)(q^{\min\{c,c_2\}}-2).$$
%\end{rmk}

\begin{lem} \label{lem:auxiliary}
If $(n,a_1,a_2,a_3)=1$, then
$$\frac{1}{\mathfrak{d}}\prod_{j=1}^3(a_j,n)=\prod_{j=1}^3(e_{j},n)=(e_1e_2e_3,n).$$
\end{lem}
\begin{proof}
For each $q\in\Upsilon(n)$, we can show
$$\sum_{j=1}^3\|(a_j,n)\|_q=\|\mathfrak{d}\|_q+\sum_{j=1}^3\|(e_j,n)\|_q,$$
in each of the following cases separately: (1) $q\nmid a_1a_2a_3$; (2) $q\mid a_j$ for exactly one $j$; (3) $q\mid a_j$ for two $j$'s.
\end{proof}

\begin{cor} \label{cor:k2}
If $(e_1,n), (e_2,n), (e_3,n)$ are pairwise coprime, then
$$\mathfrak{k}_2(n)=(e_1,n)(e_2,n)(e_3,n)-(e_1,n)-(e_2,n)-(e_3,n)+2.$$
\end{cor}

\section{Proof of Theorem \ref{thm:main}}

Suppose
$K(k/\ell)=S^2(b_1/a_1,b_2/a_2,b_3/a_3)=:M$, with $k>1$, so that $H_1(M)\cong\mathbb{Z}_{k}$.
%It is well-known that (see Section 6.3 of \cite{Sa99}) $H_1(M)$ is isomorphic to the cokernel of the endomorphism of $\mathbb{Z}^4$ defined by
%$$Q=\left(\begin{array}{cccc} a_1 & 0 & 0 & b_1 \\ 0 & a_2 & 0 & b_2 \\ 0 & 0 & a_3 & b_3 \\ 1 & 1 & 1 & 0 \\ \end{array}\right),$$
%in some natual basis. Hence
According to Section 6.3 of \cite{Sa99},
\begin{align}
k=\#H_1(M)=|b_1a_2a_3+b_2a_1a_3+b_3a_1a_2|. \label{eq:k}
\end{align}
As a consequence, $e_1,e_2,e_3\mid k$.

In the setting of Case \ref{case:1}, by (\ref{eq:equality-1'}) and Theorem \ref{thm:DW-formula},
\begin{align}
p^{h\delta_{p\mid k}}+p^h\sum\limits_{j=1}^{(n,k)-1}p^{h\omega(jn/(n,k))}=\mathfrak{k}_2(n)p^{2h}+\mathfrak{k}_1(n,p,h)p^{h}+\mathfrak{k}_0(p,h).  \label{eq:equality-1}
\end{align}
We refer to this equation as (\ref{eq:equality-1})-$(p;n,h)$.

\begin{lem}\label{lem:prepare}
Let $\alpha\in\mathbb{Z}[t]$, $d\in\mathbb{Z}_{\ge 2}$. For each $p\in\Upsilon(\Omega_d(\alpha))-\Upsilon(\alpha(1))$, there exists a divisor $d_1$ of $d$ such that $p\nmid d_1$ and
$p\mid\Omega_{d_1}(\alpha)$.
\end{lem}
\begin{proof}
Write $d=p^rd'$ with $r\in\mathbb{Z}_{\ge 0}$ and $p\nmid d'$. Suppose $\alpha(t)=a\prod_{i=1}^b(t-\lambda_i)$ over $\overline{\mathbb{F}_p}$. Then
$$\widetilde{\Omega}_{d'}(\alpha)=a^{d'}\prod\limits_{i=1}^b(\lambda_i^{d'}-1)=a^d\prod\limits_{i=1}^b(\lambda_i^d-1)=\widetilde{\Omega}_d(\alpha)=0\ \ \text{in\ \ }\mathbb{F}_p,$$
i.e. $p\mid\widetilde{\Omega}_{d'}(\alpha)$. Hence $p\mid\Omega_{d_1}(\alpha)$ for some $d_1\mid d$.
\end{proof}

\begin{rmk}
\rm As a recipe, if $p\in\Upsilon(\Omega_d(\Delta_K))$, then $\Delta_K(t)=0$ has a solution in some extension $\mathbb{F}_{p^h}$ of $\mathbb{F}_p$ that contains a primitive $d_1$-th root of unity $\xi_{d_1}$,
where $1<d_1\mid d$ and $p\nmid d_1$. This means $\Delta_K(\xi_{d_1}^v)=0$, i.e. $\omega(v)>0$, for some $v$.
\end{rmk}

When $g\mid k$, in the setting of Case \ref{case:2}, by (\ref{eq:equality-2'}) and Theorem \ref{thm:DW-formula},
\begin{align}
(p^m,k)^{o(g)}+p^{mo(g)}\sum\limits_{v=1}^{g-1}p^{o(g)\varpi_m(v)}=\mathfrak{k}_2(g)p^{2mo(g)}+\mathfrak{k}_1(g,p^m,o(g))p^{mo(g)}+\mathfrak{k}_0(p^m,o(g)).  \label{eq:equality-2}
\end{align}
We refer to this equation as (\ref{eq:equality-2})-$(p^m;g)$.

\begin{lem} \label{lem:3-coprime}
$(a_1,a_2,a_3)=1$, so that $e_1,e_2,e_3$ are pairwise coprime.
\end{lem}
\begin{proof}
Assume there exists a prime $q$ dividing $a_1,a_2,a_3$ simultaneously. Take a prime $p$ such that $p\nmid a_1a_2a_3k$ and $p>q^2$. Then (\ref{eq:equality-1})-$(p;q,h)$ becomes
$$\sum\limits_{v=1}^{q-1}p^{h\omega(v)}=(q-1)(q-2)p^h+3(q-1).$$
Clearly, $\omega(v)\in\{0,1\}$ for each $v$. Putting
$$\mathfrak{b}=\#\{v\colon 1\le v\le q-1,\ \omega(v)=1\},$$
we obtain $\mathfrak{b}p^h+q-1-\mathfrak{b}=(q-1)(q-2)p^h+3(q-1)$. This is impossible.
\end{proof}

\begin{lem} \label{lem:key}
If $e_j>1$, then $\Delta_K(t)\equiv a_j\alpha(t)\pmod{[e_j]_t}$ for some $\alpha\in\mathbb{Z}[t]$ with $\deg\alpha<e_j-1$
and $\widetilde{\Omega}_{e_j}(\alpha)=1$.
\end{lem}
\begin{proof}
Use division with remainder to write
\begin{align}
\Delta_K(t)=a'\alpha(t)+[e_j]_t\beta(t) \label{eq:remainder}
\end{align}
for some $a'\in\mathbb{Z}_{\ge 0}$, $\alpha(t),\beta(t)\in\mathbb{Z}[t]$ such that $\deg\alpha<e_j-1$, and the coefficients of $\alpha$ are coprime.

Let $p\in\Upsilon(a_j)$. By Lemma \ref{lem:3-coprime}, $p\nmid e_j$. Obviously, $p\mid a_{j-1}a_{j+1}$ if and only if $p\mid k$. For each $1<e'\mid e_j$, by (\ref{eq:equality-2})-$(p;e')$, (noting that $\mathfrak{k}_2=0$),
$$\sum_{v=1}^{e'-1}p^{o(e')\varpi_1(v)}=(e'-1)p^{o(e')}.$$
Since $p^{2o(e')}>e'p^{o(e')}$, we have $\varpi_1(v)=1$ for each $v$. In the notation introduced in Case \ref{case:2}, $b_0=1$, so that $\varpi_m(v)\le m$.
Then (\ref{eq:equality-2})-$(p^m;e')$ forces $\varpi_m(v)=m$ for each $v\in\{1,\ldots,e'-1\}$; in particular, $\Delta_K(t)=0\in\mathbb{Z}_{p^m}[t]/(\psi_{e',s}(t))$, which implies $\psi_{e',s}\mid a'\alpha$ in $\mathbb{Z}_{p^m}[t]$. Since the $\psi_{e',s}$ (for $e'\mid e$ and $s\in\Xi_p(e')$) are pairwise coprime and their product is $[e_j]_t$, we have $p^m\mid a'$.

In the other direction, whenever a prime power $p^m\mid a'$, we have $\Delta_K(t)\in(p^m,[e_j]_t)$, so that, for an arbitrary $s\in\Xi_p(e_j)$, $\Delta_K(t^v)=0\in\mathbb{Z}_{p^m}[t]/(\psi_{e_j,s})$ for each $1\le v<e_j$, then it follows from (\ref{eq:equality-2})-$(p^m;e_j)$ immediately that $p^m\mid a_j$. Therefore, $a'=a_j$.

Finally, we show $\widetilde{\Omega}_{e_j}(\alpha)=1$. By (\ref{eq:remainder}), $p\nmid\alpha(1)$. Assume on the contrary that $\Omega_d(\alpha)>1$ for some $1<d\mid e_j$. Take $p\in\Upsilon(\Omega_d(\alpha))$. We can just assume $p\nmid d$: otherwise by Lemma \ref{lem:prepare} there exists $1<d_1\mid d$ with $p\nmid d_1$ and $p\mid\Omega_{d_1}(\alpha)$. Let $m=\|a_j\|_p\ge 0$. Then $\psi_{d,s}\mid \alpha$ in $\mathbb{Z}_p[t]$ for some $s\in\Xi_p(d)$. Take $\psi\in\mathbb{Z}[t]$ which is sent by the quotient map $\mathbb{Z}[t]\twoheadrightarrow\mathbb{Z}_{p^{m+1}}[t]$ to the Hensel lift of $\psi_{d,s}$,
Writing $\alpha=p\tilde{\alpha}+\psi\eta$ in $\mathbb{Z}[t]$, we obtain from (\ref{eq:remainder}) that
$$\Delta_K(t)=pa_j\tilde{\alpha}(t)+a_j\psi(t)\eta(t)+[e_j]_t\beta(t)\in(p^{m+1},\psi(t)),$$
due to $\psi(t)\mid[d]_t$ in $\mathbb{Z}_{p^{m+1}}[t]$ (and of course $[d]_t\mid[e_j]_t$).
Now in (\ref{eq:equality-2})-$(p^{m+1};d)$, ${\rm LHS}>p^{(2m+2)o(d)}$, while ${\rm RHS}<dp^{(2m+1)o(d)}$. This is a contradiction.
\end{proof}

%\begin{rmk}
%\rm In the proof, a property of the Alexander module of $K$ is revealed.
%\end{rmk}

\begin{cor} \label{cor:determine-a}
$\Omega_d(\Delta_K)=a_j^{\phi(d)}$ for all $d\mid e_j$.
\end{cor}

\begin{lem} \label{lem:iv}
If $\widetilde{\Omega}_k(\Delta_K)=0$, then
\begin{enumerate}
  \item[\rm(a)] $e_j>1$ for at least two $j$'s;
  \item[\rm(b)] $$\frac{[e_1e_2e_3]_t}{[e_1]_t[e_2]_t[e_3]_t}\ \Big|\ \Delta_K(t);$$
%Consequently, $\Omega_{e_1e_2e_3}(\Delta_K)=0$, and $\Omega_{e_je_{j'}}(\Delta_K)=0$ if $e_j,e_{j'}>1$.
  \item[\rm(c)] if $\Omega_d(\Delta_K)=0$, then $d\mid e_1e_2e_3$.
\end{enumerate}
\end{lem}

\begin{proof}
Suppose $\Omega_{d'}(\Delta_K)=0$ for some $d'\mid k$. %Then $\Psi_d\mid \Delta_K$.
Let $p$ be an arbitrary prime with $p\nmid a_1a_2a_3k$.

(a) In (\ref{eq:equality-1})-$(p;d',h)$, $\mathfrak{k}_1$ is bounded; since $p$ can be arbitrarily large and ${\rm LHS}>p^{2h}$, we have $\mathfrak{k}_2(d')>0$. By Corollary \ref{cor:k2},
$(d',e_j)>1$ for at least two $j$'s.

(b) Take $n=e_1e_2e_3$. Then by Corollary \ref{cor:k2}, $\mathfrak{k}_2=n-e_1-e_2-e_3+2$, and (\ref{eq:equality-1})-$(p;n,h)$ forces
$\#\{i\colon 1\le i<n,\ \Delta_K(\xi_n^i)=0\}=\mathfrak{k}_2$.
By Lemma \ref{lem:key},
$$\prod_{s=1}^{e_j-1}\Delta_K(\xi_n^{se_{j-1}e_{j+1}})=a_j^{e_j-1}\ne 0\in\mathbb{F}_p$$
when $e_j>1$. Hence $\Delta_K(\xi_n^i)=0$ for all $i\in\{1,\ldots,n-1\}$ such that $ie_1,ie_2,ie_3\not\equiv 0\pmod{n}$.
Consequently,
$$\frac{[e_1e_2e_3]_t}{[e_1]_t[e_2]_t[e_3]_t}\ \Big|\ \Delta_K(t)$$
in $\mathbb{Z}_p[t]$. Since $p$ can be arbitrarily large, this actually holds in $\mathbb{Z}[t]$.

(c) Suppose $\Omega_d(\Delta_K)=0$, with $d\mid k$. Take $n=d$, and let $d_0=(e_1e_2e_3,d)=(e_1,d)(e_2,d)(e_3,d)$. Then
$$\mathfrak{k}_2=d_0-(e_1,d)-(e_2,d)-(e_3,d)+2.$$
By (\ref{eq:equality-1})-$(p;d,h)$, $\#\{i\colon 1\le i<d,\ \Delta_K(\xi_d^i)=0\}=\mathfrak{k}_2$.
From (b) we know that $\Delta_K(\xi_d^{sd/d_0})=0$ for $s\in\{1,\ldots,d_0-1\}$ which are not multiples of $d_0/(e_j,d)$ for any $j$; the number of such $s$' is exactly $\mathfrak{k}_2$. Since also $\Delta_K(\xi_d)=0$, we must have $d=d_0$, i.e. $d\mid e_1e_2e_3$.
\end{proof}

\begin{lem}
If $\widetilde{\Omega}_k(\Delta_K)\ne 0$, then $\mathfrak{k}_2(n)=0$ for any $n$, and $e_{j}>1$ for at most one $j$.
\end{lem}
\begin{proof}
Take a prime number $p$ with $p\nmid a_1a_2a_3\widetilde{\Omega}_k(\Delta_K)$.
Then $kv\equiv 0\pmod{n}$ and $\Delta_K(\xi_n^v)=0\in\mathbb{F}_{p^h}$ never hold simultaneously. So the LHS of (\ref{eq:equality-1})-$(p;n,h)$ does not exceed $np^h<p^{2h}$. Thus, $\mathfrak{k}_2=0$.

Assume $e_j>1$ for two $j$'s, say $e_1,e_2>1$. Take $q_i\in\Upsilon(e_i)$, then $q_1\ne q_2$, and we would have $\mathfrak{k}_2(q_1q_2)=(q_1-1)(q_2-1)>0$.
\end{proof}

Re-arranging $a_1,a_2,a_3$ if necessary, we assume $e_1=e_2=1$. Denote $e_3$ by $e$ and $a_3$ by $a$.
%From (\ref{eq:k}) we see $e\mid k$.

For the last term in the expression of $\mathfrak{k}_1$,
\begin{align}
\mathfrak{r}(n):=\frac{\mathfrak{e}-1}{\mathfrak{d}}\prod\limits_{j=1}^3(a_j,n)<\mathfrak{e}(e,n)=n\frac{(e,n)}{\mathfrak{d}}<n; \label{eq:last-term}
\end{align}
we have applied Lemma \ref{lem:auxiliary} and used $e<[a_1,a_2,a_3]$.

\begin{lem}\label{lem:ii}
If $\widetilde{\Omega}_k(\Delta_K)=1$, then $e=1$.
\end{lem}
\begin{proof}
Assume $e>1$. Take $p\in\Upsilon(a)$, then $p\nmid e$. In (\ref{eq:equality-1})-$(p;e,h)$, since $\Delta_K(\xi_e^v)\ne 0$ for all $v$, we have ${\rm LHS}\le e\cdot p^h\le(p^h-1)p^h$, but ${\rm RHS}>p^{2h}$, as $\mathfrak{k}_1\ge(e-1)p^h$.
This is a contradiction.
\end{proof}

\begin{lem} \label{lem:v}
Suppose $\widetilde{\Omega}_k(\Delta_K)>1$. Then
\begin{enumerate}
  \item[\rm(a)] $e>1$;
  \item[\rm(b)] if $d\mid k$ and $\Omega_d(\Delta_K)>1$, then $d\mid e$.
\end{enumerate}
\end{lem}

\begin{proof}
(a) Assume $e=1$. Since $\widetilde{\Omega}_k(\Delta_K)>1$, there exists $d\mid k$ with $\Omega_d(\Delta_K)>1$. Take $p\in\Upsilon(\Omega_d(\Delta_K))$; by Lemma \ref{lem:prepare} we may just assume $p\nmid d$.
In (\ref{eq:equality-1})-$(p;d,h)$, since $\mathfrak{k}_1=\mathfrak{r}(d)<d<p^h$ and $\mathfrak{k}_0\le p^h$,
we have ${\rm RHS}<p^{2h}$. But ${\rm LHS}\ge p^{2h}$. This is a contradiction.

(b) Take $p\in\Upsilon(\Omega_d(\Delta_K))$. By Lemma \ref{lem:prepare} there exists $d_1\mid d$ with $p\nmid d_1$ and $p\mid\Omega_{d_1}(\Delta_K)$.
By (\ref{eq:equality-1})-$(p;d_1,h)$,
$$p^{h\delta_{p\mid k}}+p^h\sum\limits_{v=1}^{d_1-1}p^{h\omega(v)}=\big(((e,d_1)-1)(a,p)^{h}+\mathfrak{r}(d_1)\big)p^{h}+\mathfrak{k}_0.$$
Consequently, $p\mid a$ and $(e,d_1)>1$. By (\ref{eq:k}), $p\nmid k$ so that actually $p\nmid d$. Replacing $d_1$ by $d$ in the above equation, and noting that $\mathfrak{k}_0=1$, we obtain
\begin{align}
\sum_{v=1}^{d-1}p^{h\omega(v)}=((e,d)-1)p^{h}+\mathfrak{r}(d).  \label{eq:consequece}
\end{align}
If $d\nmid e$, then $(e,d)<d$. By Lemma \ref{lem:key}, $[e]_t\mid\Delta_K(t)$ in $\mathbb{Z}_p[t]$, so already
$$\Delta_K(\xi_d^{sd/(e,d)})=0\in\mathbb{F}_{p^h}, \qquad  s=1,\ldots,(e,d)-1.$$
Now that $\Delta_K(\xi_d)=0$, we would have $\sum_{v=1}^{d-1}p^{h\omega(v)}>(e,d)p^{h}$. This contradicts (\ref{eq:consequece}), as $\mathfrak{r}(d)<d<p^{o(d)}$. Thus $d\mid e$.
\end{proof}

The proofs for Theorem \ref{thm:main} (i)--(v) are completed in proving Lemma \ref{lem:3-coprime}, \ref{lem:ii}, \ref{lem:key}, \ref{lem:iv}, \ref{lem:v}, respectively.

\section{Applications}

For $\alpha\in\mathbb{Z}[t]$, let $\Gamma(\alpha)=\{d\in\mathbb{Z}_{\ge 2}\colon \Omega_d(\alpha)=1\}$.

It seems to be a difficult problem to determine $\Gamma$. We independently propose
\begin{prob}
Given $\alpha\in\mathbb{Z}[t]$, how to determine $\Gamma(\alpha)$?

If $\alpha$ has a root that is not a root of unity, is $\Gamma(\alpha)$ always finite?
%Given $d$, how to find all $\alpha\in\mathbb{Z}[t]$ with $\Omega_d(\alpha)=1$?
\end{prob}

In this section, we always assume that none of the roots of $\Delta_K$ is a root of unity, so $\widetilde{\Omega}_k(\Delta_K)\ne 0$.
For convenience, call a small Seifert fibred surgery a {\it good surgery} if $\widetilde{\Omega}_k(\Delta_K)>1$. Continue to use the notations $e,a$ in the previous section.

Starting from $\Delta_K$, one can find all possible good surgeries %triples $\{k,e,a\}$ with $\widetilde{\Omega}_k(\Delta_K)>1$
through the following steps.
\begin{enumerate}
  \item Determine $\Gamma(\Delta_K)$.
  \item Find all $e\in\{2,\ldots,\deg\Delta_K+1\}$ such that $\Delta_K(t)\equiv a\alpha(t)\pmod{[e]_t}$ for some $a\in\mathbb{Z}_{\ge 2}$, $\alpha\in\mathbb{Z}[t]$ such that $\deg\alpha\le e-2$ and $\widetilde{\Omega}_e(\alpha)=1$; note that the coefficients of $\alpha$ are necessarily coprime.
  \item For each pair $e,a$, find $k$ with $e\mid k$ and satisfies: $d\in\Gamma(\Delta_K)$ for all divisors $d$ of $k$ with $d\nmid e$.
\end{enumerate}

%\begin{prop}
%If $\Delta_K$ is nontrivial, and $\Gamma(\Delta_K)$ is finite, then there are finitely many small Seifert fibred surgeries.
%\end{prop}

\begin{rmk} \label{rmk:convenient}
\rm (i) For convenient computation, we can replace the condition $\Delta_K(t)\equiv a\alpha(t)\pmod{[e]_t}$ by $(t-1)\Delta_K(t)\equiv a(t-1)\alpha(t)\pmod{t^e-1}$.

(ii) Note that for each $d\mid e$, let $\alpha_d(t)$ denote the remainder of dividing $\alpha(t)$ by $[d]_t$, then $\widetilde{\Omega}_d(\alpha_d)=1$.
This will be helpful.
\end{rmk}

Abbreviate $\Omega_n(\Delta_K)$, $\widetilde{\Omega}_n(\Delta_K)$, $\Gamma(\Delta_K)$ to $\Omega_n$, $\widetilde{\Omega}_n$, $\Gamma$, respectively.

\begin{exmp}
\rm Let $\Delta_K(t)=rt^2+(1-2r)t+r$. It is the Alexander polynomial of the twist knot $K_{2r,2}$, and also (referred to \cite{Ch21} Example 4.3) that of the odd classical pretzel knots $P(2k_1+1,2k_2+1,2k_3+1)$ with
$$1+k_1+k_2+k_3+k_1k_2+k_1k_3+k_2k_3=r.$$

We assume $r<0$ or $r\ge 3$. Now
$$\widetilde{\Omega}_s=|r^s(\lambda^s-1)(1-\lambda^{-s})|, \qquad \lambda=1+\frac{\sqrt{1-4r}-1}{2r}.$$
When $r<0$, $\widetilde{\Omega}_s$ is monotonically increasing, so $\Omega_{q^c}=\widetilde{\Omega}_{q^c}/\widetilde{\Omega}_{q^{c-1}}>1$
for each prime power $q^c$.
%$$\Omega_{q^c}=\frac{\widetilde{\Omega}_{q^c}}{\widetilde{\Omega}_{q^{c-1}}}>1.$$
Now suppose $r\ge 3$. Note that
\begin{align}
\widetilde{\Omega}_s\equiv |(2s^2-3s)r^2-2rs+1|\pmod{r^3}.  \label{eq:exmp-1}
\end{align}
If $q$ is a prime such that $\Omega_{q}=1$, then from $\Omega_q=\widetilde{\Omega}_q\equiv |2rq-1|\pmod{r^2}$ we would deduce $r\mid 2q$, and then negate (\ref{eq:exmp-1}). Hence $\Omega_{q}>1$ for each prime $q$.

From
\begin{align*}
\Delta_K(t)&\equiv 4r-1\pmod{[2]_t} \\
&\equiv (1-3r)t\pmod{[3]_t}
\end{align*}
we see that either $e=2, a=|4r-1|$, or $e=3, a=|3r-1|$.
Furthermore,
$$\Omega_4=(2r-1)^2>1, \qquad \Omega_9=|3r^3-9r^2+6r-1|>1.$$
By Lemma \ref{lem:v} (b), $k=e$.
This can be regarded as the algebraic analogue of a part of Corollary 1.2 in \cite{IM16}.

The result of Theorem 1.4 in \cite{Ka07} is quickly recovered by setting $r=-1$.
\end{exmp}

\begin{exmp}
\rm Suppose $\Delta_{K}(t)=2t^4-4t^3+5t^2-4t+2$, which is the Alexander polynomial of the knot $7_5$.
We have
\begin{align*}
\Delta_{K}(t)&\equiv 17\pmod{[2]_t} \\
&\equiv 9t^2+8\pmod{[4]_t}  \\
&\equiv 7t^2\pmod{[3]_t}  \\
&\equiv 3(-2t^3+t^2-2t)\pmod{[5]_t}.
\end{align*}
Hence the pairs $(e,a)$ probably occurring in a good surgery is $(2,17), (3,7), (5,3)$.

By computation, we find $\Omega_8>1$, $\Omega_9=1$, $\Omega_{27}>1$, $\Omega_{25}>1$. Thus, when $e=2$, it is possible for $4\mid k$, but never $8\mid k$; when $e=3$, it is possible for $9\mid k$, but never $27\mid k$; when $e=5$, it is impossible for $25\mid k$. By Lemma \ref{lem:v} (b),
$k$ cannot be a multiple of any of $6,10,15$.
To completely find all possible $k$'s, it remains to determine $\Gamma$.
\end{exmp}

%\begin{exmp}
%\rm Suppose $\Delta_K(t)=t^4-t^2+1=\Psi_{12}(t)$.

%Among $e\in\{2,3,4,5\}$, only $e=3$ is valid: $\Delta_K(t)\equiv 2(t+1)\pmod{[3]_t}$. Note that for all $c\ge 2$,
%$$\widetilde{\Omega}_{3^c}=(\zeta_{12}^{3^c}-1)(\zeta_{12}^{-3^c}-1)(\zeta_{12}^{5\cdot 3^c}-1)(\zeta_{12}^{-5\cdot 3^c}-1)=4=\widetilde{\Omega}_3,$$
%Hence $\Omega_{3^c}=1$. On the other hand, $\Omega_6=4$. Hence $6\nmid k$.

%For any $d$ with $(d,12)=1$,
%$$\widetilde{\Omega}_d=\Big|\prod\limits_{j\in\mathbb{Z}_{12}^\times}(\zeta_{12}^{dj}-1)\Big|=\Big|\prod\limits_{j\in\mathbb{Z}_{12}^\times}(\zeta_{12}^{j}-1)\Big|=|\Delta_K(1)|=1.$$
%In particular, $\Omega_q=1$ for any prime $q>3$.

%Therefore, if $K$ admits a small Seifert fibred surgery with $k>1$ and $e_j>1$ for some $j$, then such $j$ is unique, in which case $a_j=2$, $e_j=3$, and $k=3^c$ for some $c\ge 1$.
%\end{exmp}

\begin{exmp}
\rm Suppose
$$\Delta_K(t)=t^{10}-t^9+t^7-t^6+t^5-t^4+t^3-t+1,$$
which is the Alexander polynomial of $P(-2,3,7)$, the $(-2,3,7)$-pretzel knot.
Direct computations show $\Omega_2=\Omega_3=1$ and
\begin{align*}
\Delta_K(t)&\equiv 3t^3+3t^2+5\pmod{[5]_t}  \\
&\equiv 2t^5+3t^3+3\pmod{[7]_t} \\
&\equiv -2t^9-t^8-2t^6-2t^4-t^2-2t\pmod{[11]_t}.
\end{align*}
By Remark \ref{rmk:convenient} (ii), we need not to check the remainders of dividing $\Delta_K(t)$ by $[10]_t$.
Thus, $K$ does not admit a good surgery.

Actually, by \cite{Me14} Theorem 1.1, $P(-2,3,7)$ admits no small Seifert fibred surgery.

%Generally, applying the method developed in \cite{Ch21}, we can compute
%\begin{align}
%\Delta_{P(-2,3,2r+1)}(t)=(t^4-t^2+t)[2r+1]_{-t}+(t-1)(t-1-t^2).  \label{eq:even-pretzel}
%\end{align}
%Therefore, when $\Delta_K(t)$ equals the RHS of (\ref{eq:even-pretzel}),
\end{exmp}

\begin{exmp}
\rm Let $K$ be the Montesinos knot $M(-2,5/2,2r+1)$ (in the notation of \cite{Ch21}, but different from the notation usually used in the literature).
Applying the method of \cite{Ch21}, we compute
$$\Delta_K=-t^{2r+2}+2t^{2r+1}-t^3\frac{t^{2r-3}+1}{t+1}+2t-1.$$

Consider the two cases $r=1$ and $r=2$. We only write down the relevant part of computation.

When $r=1$,
\begin{align*}
\Delta_K&=-t^4+2t^3-t^2+2t-1 \\
&\equiv -7\pmod{[2]_t} \\
&\equiv 2(t+1)\pmod{[3]_t} \\
&\equiv -3t^2-4\pmod{[4]_t} \\
&\equiv 3(t^2+t)\pmod{[5]_t}.
\end{align*}
This is consistent with $K(3)=S^2(1/2,-1/3,-2/15)$, $K(4)=S^2(1/2,-1/6,-2/7)$, $K(5)=S^2(-1/3,-1/5,3/5)$, as in \cite{Me14} Theorem 1.4.

When $r=2$,
\begin{align*}
\Delta_K=-t^6+2t^5-t^3+2t-1\equiv 5t\pmod{[4]_t}.
\end{align*}
As easily computed, $\Omega_8=1$. By \cite{Me14} Theorem 1.4, $K(8)=S^2(-1/4,3/4,-2/5)$. This example is given for showing that $(e,k/e)>1$ really may occur.
\end{exmp}

\section{Proof of Theorem \ref{thm:DW-formula}}

A formula for the Dijkgraaf-Witten invariant of a general orientable Seifert 3-manifold was derived in \cite{Ch12}. For each conjugacy class of a finite group $G$, choose a representative $x$; for each isomorphism class of irreducible representations of ${\rm Cen}(x)$, choose a representative $\rho$, with character $\chi_\rho$. Let $\Lambda$ denote the set of all such pairs.
For $(x,\rho)\in\Lambda$, and $a,b\in\mathbb{Z}$ with $(a,b)=1$, put
\begin{align*}
\eta_{a,b}(x,\rho)=\frac{1}{\dim\rho}\cdot\sum\limits_{z\in{\rm Cen}(x)\atop z^a=x}\chi_{\rho}(z^{-b}).
\end{align*}
By (32) of \cite{Ch12},
\begin{align*}
\#\hom\big(\pi_1\big(S^2(b_1/a_1,b_2/a_2,b_3/a_3)\big),G\big)=\# G\cdot\sum\limits_{(x,\rho)\in\Lambda}T(x,\rho),
\end{align*}
where
\begin{align*}
T(x,\rho)=\left(\frac{\dim\rho}{\#{\rm Cen}(x)}\right)^{2}\cdot\prod\limits_{j=1}^{3}\eta_{a_j,b_j}(x,\rho).
\end{align*}

Recall the construction in Section 8.2 of \cite{GTM42}.
Suppose $G=A\rtimes H$ with $A$ abelian. Let $H$ act on $A^\vee:=\hom(A,\mathbb{C}^\times)$ by $(h\mu)(a)=\mu(h^{-1}.a)$, for $h\in H, \mu\in A^\vee$, and let $\{\mu_i\colon i\in A^\vee/H\}$ be a system of representatives for the orbits.

For each $i$, let
$$H_i=\{h\in H\colon h\mu_i=\mu_i\}=\{h\in H\colon \mu_i(h^{-1}.a)=\mu_i(a)\ \text{for\ all\ }a\in A\}.$$
Extend $\mu_i$ to $G_i:=A\cdot H_i$ by $\mu_i(ah)=\mu_i(a)$.
Let $\rho$ be an irreducible representation of $H_i$; it gives rise to a representation $\tilde{\rho}$ by post-composing with $G_i\twoheadrightarrow H_i$. Let $\theta_{i,\rho}$ denote the representation of $G$ induced from $\mu_i\otimes\tilde{\rho}$.

By Proposition 25 of \cite{GTM42}, each irreducible representation of $G$ is isomorphic to a unique one of the $\theta_{i,\rho}$'s.

\bigskip

Now let $A=\mathbb{Z}_{p^m}^h$, and $G=A\rtimes_f\mathbb{Z}_{n}$. There is a canonical bilinear pairing $A\times A\to \mathbb{Z}_{p^m}$ which we denote by $w\cdot w'\in\mathbb{Z}_{p^m}$, for $w,w'\in A$.

As easily seen, the nontrivial conjugacy classes of $G$ are
${\rm Con}(\beta^v)$, $v=1,\ldots,n-1$ and ${\rm Con}(u)$, $u\in O$,
where $O$ is a system of representatives for the orbits of $A-\{0\}$ under the action of $f$.
The centralizers of the representatives are: ${\rm Cen}(u)=A$ for $u\in O$; ${\rm Cen}(\beta^v)=\langle\beta\rangle\cong\mathbb{Z}_n$ for $1\le v\le n-1.$

For $w\in A$, define $w^\vee\in A^\vee$ by sending $u$ to $\zeta_{p^m}^{w\cdot u}$.
For $t\in\mathbb{Z}_n$, define $t^\vee\in\mathbb{Z}_n^\vee$ by sending $\beta$ to $\zeta_n^{t}$.

The action of $\langle\beta\rangle=\mathbb{Z}_n$ on $\{w^\vee\colon w\in A\}$ is given $\beta^vw^\vee=(\phi^{-v}w)^\vee.$
A system of representatives for the orbits can be taken as $\{1\}\cup\{w^\vee\colon w\in O\}$.
Obviously, the stabilizer of each $\mu_w$ is trivial.
So the irreducible representations of $G$ are $\varrho_t$, $t\in\mathbb{Z}_n$, and $\mu_w$, $w\in O$,
where $\varrho_t$ is the composite of $t^\vee$ with $G\twoheadrightarrow \mathbb{Z}_n$, and $\mu_w$ is the representation induced by $w^\vee$.

Then $\Lambda=\Lambda_1\sqcup\Lambda_2\sqcup\Lambda_3$, with
\begin{align*}
\Lambda_1&=\{(u,w^\vee)\colon u\in O, w\in A\}, \\
\Lambda_2&=\{(\beta^v,t^\vee)\colon 1\le v\le n-1,\ t\in\mathbb{Z}_n\}, \\
\Lambda_3&=\{(1,\varrho_t)\colon 1\le t\le n\}\sqcup\{(1,\mu_w)\colon w\in O\}.
\end{align*}
Let $d=[a_1,a_2,a_3]$. Abbreviate $\|a\|_p$ to $\|a\|$. %Let $z=\max\{\|a_1\|,\|a_2\|,\|a_3\|\}$.
\begin{itemize}
  \item For each $u\in O$ and each $w\in A$,
        $$\eta_{a,b}(u,w^\vee)=\sum\limits_{x\in A\atop ax=u}\zeta_{p^m}^{-bw\cdot x}
        =(a,p^{m})^h\delta_{\|u\|\ge\|a\|}\delta_{\|w\|\ge\|a\|}\zeta_{p^m}^{-(bw/a)\cdot u}.$$
        Hence
        $$T(u,w^\vee)=\frac{1}{p^{2mh}}\delta_{z<m}\delta_{\|u\|\ge z}\delta_{\|w\|\ge z}\zeta_{p^m}^{-(\mathfrak{n}w/d)\cdot u}
        \prod\limits_{j=1}^3(a_j,p^{m})^h.$$
        When $\|u\|,\|w\|\ge z$, writing $u=p^zu_1$, $w=p^zw_1$,
        \begin{align*}
        \sum\limits_{w\in A\atop\|w\|\ge z}\zeta_{p^m}^{-(zw/d)\cdot u}
        =\sum\limits_{w_1}\zeta_{p^{m-z}}^{-\mathfrak{n}w_1u_1p^z/d}
        =p^{h(m-z)}\delta_{\|u_1\|\ge m-z-\|\mathfrak{n}\|},
        \end{align*}
        implying
        \begin{align}
        &\sum\limits_{u\in O^t}\sum\limits_{w\in A}T(u,w^\vee) \nonumber \\
        =&\ \frac{\delta_{z<m}}{np^{h(m+z)}}\big(\#\{u\in A\colon \|u\|\ge\max\{z,m-\|\mathfrak{n}\|\}\}-1\big)\prod\limits_{j=1}^3(a_j,p^m)^h \nonumber \\
        =&\ \delta_{z<m}\frac{p^{h\min\{m-z,\|\mathfrak{n}\|\}}-1}{np^{h(m+z)}}\prod\limits_{j=1}^3(a_j,p^{m})^h.  \label{eq:proof-1}
        \end{align}
  \item For each $t$,
        $$\eta_{a,b}(\beta^v,t^\vee)=\sum\limits_{y\in\mathbb{Z}_n\atop ay=v}\zeta_n^{-tby}
        =\delta_{(a,n)\mid v}\delta_{(a,n)\mid t}\zeta_n^{-(tb/a)v},$$
        Hence
        \begin{align*}
        \sum\limits_{t=1}^nT(\beta^v,t^\vee)=\sum\limits_{t=1}^n\frac{\delta_{\mathfrak{d}\mid v}\delta_{\mathfrak{d}\mid t}\zeta_n^{-tv\mathfrak{n}/d}}{n^2}\prod\limits_{j=1}^3(a_j,n)
        =\frac{\delta_{\mathfrak{d}\mid v}\delta_{n\mid \mathfrak{n}v}}{n\mathfrak{d}}\prod\limits_{j=1}^3(a_j,n).
        \end{align*}
        Thus, using $\delta_{\mathfrak{d}\mid v}\delta_{n\mid \mathfrak{n}v}=\delta_{[\mathfrak{d},n/(n,\mathfrak{n})]\mid v}$ and $n=\mathfrak{e}[\mathfrak{d},n/(n,\mathfrak{n})]$, we obtain
        \begin{align}
        \sum\limits_{v=1}^{n-1}\sum\limits_{t=1}^nT(\beta^v,t^\vee)=\frac{\mathfrak{e}-1}{n\mathfrak{d}}\prod\limits_{j=1}^3(a_j,n).  \label{eq:proof-2}
        \end{align}
  \item Note that, if $0\ne v\in\mathbb{Z}_n$, then $(u\beta^v)^a=1$ is equivalent to $dv=0$. So
        \begin{align*}
        \eta_{a,b}(1,\varrho_t)&=\#\{u\in A\colon au=0\}+\sum\limits_{0\ne v\in\mathbb{Z}_n\atop av=0}\sum\limits_{u\in A}\zeta_n^{-btv} \\
        &=(a,p^m)^h-p^{mh}+p^{mh}(a,n)\delta_{(a,n)\mid t}.
        \end{align*}
        Consequently,
        \begin{align*}
        T(1,\varrho_t)=\frac{1}{p^{2mh}n^2}\prod\limits_{j=1}^3\big((a_j,p^m)^h-p^{mh}+p^{mh}(a_j,n)\delta_{(a_j,n)\mid t}\big),
        \end{align*}
        and then
        \begin{align}
        \sum\limits_{t=1}^nT(1,\varrho_t)=\frac{1}{p^{2mh}n}\Big(\mathfrak{k}_2p^{3mh}+\mathfrak{k}'p^{2mh}+\prod\limits_{j=1}^3(a_j,p^m)^h\Big),  \label{eq:proof-3}
        \end{align}
        with
        $$\mathfrak{k}'=\sum\limits_{j=1}^3((a_{j+1},a_{j-1},n)-1)(a_j,p^m)^h.$$
  \item Finally, noting
        $\chi_{\mu_w}(u\beta^v)=\delta_{v=0}\cdot\sum_{y=1}^n\zeta_{p^m}^{w\cdot f^{-y}u},$
        we have
        \begin{align*}
        \eta_{a,b}(1,\mu_w)
        =\frac{1}{n}\sum\limits_{u\in A\atop au=0}\sum\limits_{y=1}^n\zeta_{p^m}^{w\cdot f^{-y}u}
        =\sum\limits_{u\in A\atop au=0}\zeta_{p^m}^{w\cdot u}
        =(a,p^m)^h\delta_{\|w\|\ge\|a\|}.
        \end{align*}
        Thus $T(1,\mu_w)=p^{-2mh}\delta_{\|w\|\ge z}\prod_{j=1}^3(a_j,p^m)^h$, implying
        \begin{align}
        \sum\limits_{w\in O^t}T(1,\mu_w)=\delta_{z<m}\frac{p^{h(m-z)}-1}{p^{2mh}n}\prod\limits_{j=1}^3(a_j,p^m)^h. \label{eq:proof-4}
        \end{align}
\end{itemize}

The formula is established by summing (\ref{eq:proof-1})--(\ref{eq:proof-4}) and multiplying by $p^{mh}n$.

{\large Haimiao Chen}  \ \ \ Email: chenhm$@$math.pku.edu.cn \\
Mathematics, Beijing Technology and Business University, 11 Fucheng Road, Haidian District, 100048, Beijing, China

\end{document}